\documentclass[12pt]{article}
\usepackage{amsthm}
\usepackage{amsmath}
\usepackage{amssymb}
\usepackage{amsfonts}
\usepackage{latexsym}
\newtheorem{thm}{Theorem}[section]
\numberwithin{equation}{section}
\newcommand{\nc}{\newcommand*}
\newcommand{\rnc}{\renewcommand*}
\rnc{\baselinestretch}{1.2}
\def\reff#1{(\ref{#1})}
\def\ds{\displaystyle}
\nc{\ts}{\textstyle}

\def\half{{\frac{1}{2}}}
\nc{\ket}[1]{{\vert{#1}\rangle}}               
\nc{\bra}[1]{{\langle{#1}\vert}}               
\nc{\braket}[2]{{\langle{#1}\vert{#2}\rangle}} 
\nc{\scpr}[2]{\left({#1}\,,\,{#2}\right)}      
\nc{\ccr}[2]{\left[{#1}\,,\,{#2}\right]}       
\nc{\car}[2]{\left\{{#1}\,,\,{#2}\right\}}     
\nc{\pois}[2]{\left\{{#1},{#2} \right\}_{cl}}  
\nc{\Szi}[1]{\sum_{ #1 = 0 }^{\infty}}
\nc{\Szx}[2]{\sum_{#1 =0}^{#2}}
\nc{\Syx}[3]{\sum_{ #1 = #2 }^{#3}}
\nc{\goth}{\mathfrak}
\rnc{\bold}{\mathbf}
\rnc{\frak}{\mathfrak}
\rnc{\Bbb}{\mathbb}
\rnc{\rm}{\text}
\nc{\mbf}[1]{\mbox{\boldmath ${#1}$}}
\nc{\BR}{\mathbb R}
\nc{\BC}{\mathbb C}
\nc{\BI}{{\mathbb I}}

\def\cN{{\mathcal N}}
\nc{\gmf}[1]{\varGamma({ #1})}
\nc{\Foo}[3]{\raisebox{-3pt}{${}_1\text{\large F}_1$}%
 \left(\ts{\genfrac{}{}{0pt}{}{#1}{#2}\left.\right| #3}\right)}
\nc{\Fto}[4]{\raisebox{-3pt}{${}_2\text{\large F}_1$}%
 \left(\ts{\genfrac{}{}{0pt}{}{#1, #2}{#3}\left.\right| #4}\right)}
\nc{\Ftf}[7]{\raisebox{-3pt}{${}_2\text{\large F}_4$}%
 \left(\ts{\genfrac{}{}{0pt}{}{#1, #2}{#3, #4, #5, #6}%
 \left.\right| #7}\right)}
\begin{document}

{\begin{center} {\Huge\bf Generalized Coherent States for }\\[.5cm]
{\Huge\bf Classical Orthogonal Polynomials}
{\huge\footnote{This research was supported by RFFI grant No
00-01-00500}}
\end{center}
\vskip .3cm

\centerline{\Large\bf V.V.Borzov${}^{\left.*\right)}$,\ and
E.V.Damaskinsky${}^{\left.**\right)}$}
\begin{center}
{${}^{\left.*\right)}$\ Department of Mathematics,\\ St.Petersburg
University of Telecommunications, \\ 191065, Moika  61,
St.Petersburg, Russia\\ E-mail address: vadim@VB6384.spb.edu}
\\
{${}^{\left.**\right)}$\ University of Defense Technical
Engineering,\\Zacharievskaya 22, St.Petersburg, Russia\\ E-mail
address: evd@ED6911.spb.edu and evd@pdmi.ras.ru}
\end{center}

\vskip 0.5cm

\centerline{\bf Abstract}
\begin{quote}
For the oscillator-like systems, connected with the Laguerre, Legendre
and Chebyshev polynomials coherent states of Glauber-Barut-Girardello
type are defined. The suggested construction can be applied to
each system of orthogonal polynomials including classical ones as well as
deformed ones.
\end{quote}
\vskip 0.5cm

\section{Introduction}
We consider the oscillator-like systems, connected with the
Laguerre, Legendre and Chebyshev polynomials in the same way as
usual boson oscillator connected with the Hermite polynomials. It
is natural to call these systems as the Laguerre, Legendre and
Chebyshev (generalized) oscillators, respectively. We define
analogues of coherent states for these generalized oscillators. It
is possible to consider such definitions as a new method for
construction of coherent states. In general, this method can be
applied for a  orthogonal functions  system connected with a
Jacobi matrix.

There are four different  definitions of coherent states
\begin{enumerate}
\item Glauber-Barut-Girardello coherent states;
\item Perelomov-Gilmore coherent states;
\item Minimum uncertainty  coherent states;
\item Klauder-Gazeau coherent states.
\end{enumerate}
Here we will consider only first of these definitions.

Let $a$ and $a^{\dagger}$ are the ladder operators (annihilation
and creation, respectively) on the Fock space $\mathcal F$ with
standard orthonormal  basis $\left\{
{\ket{n}}\right\}_{n=0}^{\infty }.$

The Glauber-Barut-Girardello coherent states are defined as
eigenstates of annihilation operator
\begin{equation}
a\ket{z}=z\ket{z},\quad {z\in{\BC}}.
\label{triv1}
\end{equation}

It is well known that
\begin{equation}
\ket{z}=\exp(-\half|z|^2)\Szi{n}\frac{z^n}{\sqrt{n!}}\ket{n}=
    \exp(-\half |z|^2)\exp(za^{\dagger})\ket{0},
\label{triv2}
\end{equation}
where $\ket{0}$ is the Fock vacuum, such that $a\ket{0}=0$ and
\begin{equation}
\ket{n}=\ds\Szi{n}\frac{{a^{\dagger}}^n}{\sqrt{n!}}\ket{0},\quad
{n=0,1,...\quad  .}
\label{triv3}
\end{equation}

Let us remark that in the case of standard boson oscillator all of the
above mentioned definitions are equivalent. It means that each of
them generates the same set of coherent states. But this is not
true in the general case.

In the following we extend the above definition of coherent states
to the cases of  Laguerre, Legendre and Chebyshev (generalized)
oscillators. In doing this we use mainly the  generating
functions method and the theory of classical power moment problem for
Jacobi matrix, connected with related orthonormal polynomial systems.

\section{The generalized oscillator algebra}

Let $\mu$ denotes a positive Borel measure on the real line ${\BR}^1$
for which all moments $\mu_k$
\begin{equation}
\mu_0=1,\qquad
\mu_{k}=\int_{-\infty}^{\infty}{x^{k}\mu(dx)},\qquad k= 0,1,\dots.
\label{triv6}
\end{equation}
are finite. We also suppose that measure $\mu$ is a symmetric one,
that is
\begin{equation}
\mu_{2k+1}=0,\qquad k= 0,1,\dots  .
\label{triv7}
\end{equation}

Let us consider a system
$\left\{ {\Psi_{n}(x)}\right\}_{n=0}^{\infty}$
of polynomials defined by the recurrence relations
($n\geq 0$):
\begin{equation}
{x {\Psi_{n}(x)}}=
{b_{n}{\Psi_{n+1}(x)}}+{b_{n-1}{\Psi_{n-1}(x)}}, \qquad
\Psi_{0}(x)=1,\qquad{ b_{-1}}=0 ,
\label{triv8}
\end{equation}
where $\left\{ {b_{n}}\right\}_{n=0}^{\infty }$ is a given positive
sequence.

The following theorem was proved in \cite{bor1}
\begin{thm} \label{a}
 The polynomial system
 $\left\{ {\Psi_{n}(x)}\right\}_{n=0}^{\infty }$
is orthonormal one in the Hilbert  space
${\mathcal H}_x={\rm L}^2(R^1;{\mu)}$ if and only if the
coefficients $b_n$ and the moments $\mu_{2k}$ are connected by the
following relations
\begin{equation}
\label{triv9}
 \Szx{m}{[\half n]}\Szx{s}{[\half n]}
(-1)^{m+s}\alpha_{2m-1, n-1}\alpha_{2s-1, n-1}
\frac{\mu_{2n-2m-2s+2}}{(b_{n-1}^{\,2})!}=
b_{n-1}^{\,2}+b_{n}^{\,2},\qquad n= 0,1,\dots  .
\end{equation}
where $(b_{n}^{\,2})!=b_{0}^{\,2}b_{1}^{\,2}\cdots b_{n}^{\,2},$
the integral part of $a$ is  denoted by $[a]$,  and the
coefficients $\alpha_j$ are given by
\begin{equation}
\alpha_{2p-1, n-1}= \sum_{k_1=2p-1}^{n-1} b_{k_1}^{\,2}
\sum_{k_2=2p-3}^{k_1-2} b_{k_2}^{\,2} \cdots
\sum_{k_p=1}^{k_{p-1}-2}b_{k_p}^{\,2}.
 \label{triv10}
\end{equation}
\end{thm}
The recurrence relations (\ref{triv8}) determine the position
operator $X,$  realized as the operator of
multiplication by argument in the space ${{\mathcal H}_x}.$
We define the momentum operator $P$
and Hamiltonian $H$ by the following relations
\begin{equation}
P=K^{*}YK,\qquad H=X^{2}+P^{2}.
 \label{triv11}
\end{equation}
where $Y$ is the position operator in the (momentum) space
${\mathcal H}_y,$ and
$K: {\mathcal H}_x\longmapsto {\mathcal H}_y$
and
${K^{*}}: {\mathcal H}_y\longmapsto{\mathcal H}_x$
are the unitary integral operators  with Poisson kernels
\begin{equation}
{\mathfrak{K}(x,y;-i)}= \sum_{n=0}^{\infty }{{(-i)^{n}}\cdot
{\varphi_{n}(x)}\cdot {\varphi_{n}(y)}}. \label{triv12}
\end{equation}
and $\overline{\mathfrak{K}}$, respectively.

In what follows we consider the Hilbert space ${{\mathcal H}_x}$
as a Fock space with the basis\hfill\break
$\left\{{\ket{n}=\Psi_{n}(x)}\right\}_{n=0}^{\infty }$ and we
define the creation $a^{\dagger}$ and annihilation $a$ operators
by the formulas
\begin{equation}
a^{+}=\frac{1}{\sqrt{2}}\left( X+\imath P \right),\qquad
a=\frac{1}{\sqrt{2}}\left( X-\imath P \right).
 \label{triv13}
\end{equation}
We also define the number operator $N$
\begin{equation}
N\ket{n}={n}\ket{n},\qquad n\geq  0.
\label{triv14}
\end{equation}
and the operator-function $B(N)$ acting on the basis states by
\begin{equation}
B(N)\ket{n}=b_{n-1}^{2}\ket{n},\qquad n\geq  0,\quad b_{-1}=0.
\label{triv15}
\end{equation}

The following theorem can be proved by direct calculation.
\begin{thm}
\label{b}
 The operators $a=a^{-}, a^{+}$ and $N$ obey the
following commutation relations
\begin{equation}
\ccr{a^{-}}{a^{+}}=2\left( B(N+I)-B(N) \right),\quad
\ccr{N}{a^{\pm}}=\pm a^{\pm}. \label{triv42}
\end{equation}
Moreover if there exists real number $A$ and real function $C(n)$,
such that
\begin{equation}
b_{n}^{2}-Ab_{n-1}^{2}=C(n),\qquad n\geq  0,\quad b_{-1}=0,
\label{triv43}
\end{equation}
then ladder operators $a^{\pm}$ fulfill the relations
\begin{equation}
a^{-}a^{+}-Aa^{+}a^{-}=2C(N). \label{triv44}
\end{equation}
\end{thm}
\medskip

We call so obtained algebra $A_{\Psi}$ as generalized oscillator
algebra connected with given system of orthogonal polynomials.

{\bf Remark}. In the case of a non symmetric measure (with
$\mu_{2k+1}\neq 0$) the recurrent relations are more complicated
\begin{equation}
{x {\Psi_{n}(x)}}=
{b_{n}{\Psi_{n+1}(x)}}+{a_{n}{\Psi_{n}(x)}}+{b_{n-1}{\Psi_{n-1}(x)}},
\qquad \Psi_{0}(x)=1,\qquad{ b_{-1}}=0,
 \label{triv5a}
\end{equation}
so that the related Jacobi matrix has nonzero diagonal. In this
case one can repeat the full construction with necessary
modifications and complications (see \cite{bor1}). Such situation
arises, for example, in the cases of the Laguerre and Jacobi
polynomials.

\section{Coherent states connected with classical orthogonal polynomials}

\subsection{Coherent states of Glauber-Barut-Girardello type}

For a system of polynomials
$\left\{ {\Psi_{n}(x)}\right\}_{n=0}^{\infty }$ in the space
${{\mathcal H}_x}$
we define the coherent states of Glauber-Barut-Girardello type by
the relation
\begin{equation}
\label{triv16}
\ket{z}={\cN}(|z|^2)\Szi{n}\frac{z^n}{(\sqrt{2}b_{n-1})!}\ket{n},
\end{equation}
where the coefficients in the recurrent relations \reff{triv8} (or
\reff{triv5a} ) are denoted by $b_n$  and normalizing factor is
given by the relation
\begin{equation}
{\cN}^2=\braket{z}{z}=\Szi{n}\frac{|z|^{2n}}{(2b_{n-1}{}^2)!}\Doteq
{\exp}_{[2b_{n-1}^2]}(|z|^{2n}).
\label{triv17}
\end{equation}
Below we apply the above construction to the cases in which some
well known classical polynomials are taken as given system of
polynomials $\left\{ {\Psi_{n}(x)}\right\}_{n=0}^{\infty }$

\subsubsection{Laguerre coherent states of Glauber-Barut-Girardello type}

Let $\left\{ {\Psi_{n}(x)}\right\}_{n=0}^{\infty }$ be the orthogonal
system of Laguerre polynomials in the Hilbert space
${\mathcal H}={\rm L}^2(R^1;{\mu)}$ with the probability measure
\begin{equation}
\mu(\text{d}x)=x^{\alpha}e^{-x}\frac{\text{d}x}{\sqrt{\varGamma
(\alpha+1)}}.
\end{equation}
By definition
\begin{equation}
\Psi_{n}(x)\Doteq\sqrt{\frac{\varGamma (\alpha+n+1)}{n!\varGamma
(\alpha+1)}} \Foo{-n}{\alpha+1}{x},\qquad n\geq 0.
\label{triv18}
\end{equation}
The coefficients $a_n,b_n$ in recurrent relations {\reff{triv5a}}
 are equal to
\begin{equation}
a_n=2n+\alpha+1,\quad
b_n=-\sqrt{(n+1)(n+\alpha+1)},\qquad n\geq 0.
\label{triv19}
\end{equation}
The normalizing factor
\begin{equation}
 {\cN}^2=\left(\frac{2}{z}\right)^{\alpha} \varGamma (\alpha+1)
I_{\alpha}(\sqrt{2}|z|),
\label{triv20}
\end{equation}
  where
\begin{equation}
I_{\alpha}(z)=\left(\frac{z}{2}\right)^{\alpha}\Szi{n}
\frac{z^n}{2^nn!\varGamma (\alpha+n+1)}
\label{triv21}
\end{equation}
is the Bessel function of the first kind. According to
{\reff{triv16}}, we have
\begin{equation}
\ket{z;\alpha}=\frac{1}{(\sqrt{2}x)^{\alpha}
I_{\alpha}(\sqrt{2}z)} e^{z/\sqrt{2}}
I_{\alpha}(2^{3/4}\sqrt{xz}).
\label{triv22}
\end{equation}
The overlap of two coherent states is equal to
\begin{equation}
\braket{z_1;\alpha}{z_2;\alpha}=
\frac{I_{\alpha}(2\sqrt{\overline{z_1} z_2)}}
{\sqrt{I_{\alpha}(\sqrt{2}|z_1|)I_{\alpha}(\sqrt{2}|z_2|)}}
\label{triv23}
\end{equation}
One can check the validity of unity decomposition
\begin{equation}
 \int
\ket{z;\alpha}\bra{z;\alpha}\text{d}\nu(z,\alpha)=\BI,
\label{triv24}
\end{equation}
with the measure
\begin{equation}
\text{d}\nu(z,\alpha)=\frac{\sqrt{2}}{\pi}K_{\alpha}(\sqrt{2}|z_2|)
I_{\alpha}(\sqrt{2}|z_2|),
 \label{triv25}
\end{equation}
where $K_{\alpha}$ is a modified  Bessel function of the second
kind.

Let us note that this result coincide with results of other
authors \cite{kida} obtained by different methods.

\subsubsection{Legendre coherent states of Glauber-Barut-Girardello type}

We consider orthogonal system $\left\{
{\Psi_{n}(x)}\right\}_{n=0}^{\infty }$ of Legendre polynomials
\begin{equation}
 \Psi_{n}(x)\Doteq\sqrt{2n+1} \Fto{-n}{n+1}{1}{\frac{1-x}{2}},
\label{triv26}
\end{equation}
in the Hilbert space
${\mathcal H}={\rm L}^2([-1,1];\half\text{d}x).$

The coefficients $b_n$ in the recurrent relations {\reff{triv8}}
are equal to
\begin{equation}
b_n\Doteq \sqrt{\frac{(n+1)^2}{(2n+1)(2n+3)}},\qquad n\geq 0.
\label{triv27}
\end{equation}
Then
\begin{equation}
  (2b_{n-1}^2)!=\frac{n!(1)_n}{(\half)_n(\frac32)_n},\qquad n\geq
  1,
 \label{triv28}
\end{equation}
where Pochhammer symbol  $(a+1)_n$ is defined by
\begin{equation}
(a)_0=1,\quad
(a)_n=a(a+1)\cdots (a+n-1)=\frac{\gmf{a+n}}{\gmf{a}}, \quad
n=1,2,\ldots .
 \label{triv29}
\end{equation}
The normalizing factor is equal to
\begin{equation}
 {\cN}^2=\Szi{n}
\frac{(1/2)_n(3/2)_n}{n!(1)_n}(2|z|^2)^n =
\Fto{1/2}{3/2}{1}{2|z|^2}.
\label{triv30}
\end{equation}
Note that the radius of convergence for series in the above
relation is equal to $1/\sqrt{2}.$

For coherent state of the Legendre oscillator one obtains
\begin{equation}
\ket{z}=\frac{\Fto{3/4}{5/4}{1}{\frac{(x^2-1)4z^2}{(1-2xz)^2}}}
{\sqrt{\Fto{1/2}{3/2}{1}{2|z|^2}}} (1-2xz)^{-3/2}.
\label{triv31}
\end{equation}
The overlap of two coherent states is equal to
\begin{equation}
\braket{z_1}{z_2}= \frac{\Fto{1/2}{3/2}{1}{2\overline{z_1} z_2}}
{\sqrt{\Fto{1/2}{3/2}{1}{2|z_1|^2}\Fto{1/2}{3/2}{1}{2|z_2|^2}}}.
 \label{triv32}
\end{equation}
Finally, we have
\begin{equation}
 \int_{{\BC}_{1/\sqrt{2}}}
\ket{z}\bra{z}\text{d}\nu(z)=\BI,
 \label{triv33}
\end{equation}
with the measure
\begin{equation}
\text{d}\nu(z)=
\frac{(4|z|^2-5)P_{1/2}(|z|^2-1)-3P_{3/2}(|z|^2-1)}{2(|z|^2-2)}
\text{d}(\text{Re}z)\text{d}(\text{Im}z),\quad 0<|z|< 1/\sqrt{2},
 \label{triv34}
\end{equation}
where by $P_{\alpha}(x)$ we denote the  Legendre function.

\subsubsection{Chebyshev coherent states of Glauber-Barut-Girardello type}

The coherent states of Glauber-Barut-Girardello type for the
Chebyshev oscillator in the space
\begin{equation}
 {\mathcal H}={\rm
L}^2([-1,1];\frac{\text{d}x}{\pi\sqrt{1-x^2}})
 \label{triv35}
\end{equation}
looks as (see (\cite{bd1}))

\begin{equation}
\ket{z}_1=\frac{\sqrt{2}}{1-2|z|^2}
\frac{1-\sqrt{2}zx}{1-2\sqrt{2}zx+2z^2}, \qquad |z|<1/\sqrt{2}.
 \label{triv36}
\end{equation}

The related holomorphic representation in the Bargmann type space
consists from functions analytical in the disk $|z|<1/\sqrt{2}.$
The constructed family of coherent states possess all standard
properties. The decomposition of unity defines measure represented
by $\delta$-function on the boundary of the disk $|z|<1/\sqrt{2}.$

\end{document}